%% file: BI+MS.tex
 \documentclass[graybox]{svmult}
 \usepackage{algorithm}
 \usepackage{romannum}

 \usepackage{type1cm}        

 \newcommand{\xx}{\boldsymbol{x}}

 \newcommand{\C}{\mathcal{C}}

 \newcommand{\N}{\mathcal{N}}
 
 \usepackage{makeidx}         
 \usepackage{graphicx}        
 \usepackage{multicol}        
 \usepackage[bottom]{footmisc}
 \newcommand{\RNum}[1]{\uppercase\expandafter{\romannumeral #1\relax}}
 
 \DeclareRobustCommand{\BcalF}{{\boldsymbol{\mathcal F}}}
 \newcommand{\blue}[1]{\textcolor{black}{#1}}
 \usepackage{newtxtext}       %
 \usepackage{newtxmath}       
 
 \makeindex             

 \usepackage{xcolor}

 \newcommand{\Blambda    }{{\boldsymbol{\lambda}}}

 \newcommand{\calB}{\mathcal{B}}
 \newcommand{\calT}{\mathcal{T}}

 \newcommand{\calR}{\mathcal{R}}

 \newcommand{\calH}{\mathcal{H}}
 
 \newcommand{\bs}[1]{{\mathbf{#1}}} 
 \newcommand{\Bu}{\bs{u}}
 \newcommand{\Bn}{\bs{n}}

 \newcommand{\Bx}{\bs{x}}

 \newcommand{\BI}{\bs{I}}

 \newcommand{\Bvarphi}{{\boldsymbol{\varphi}}}
 \newcommand{\BX}{\bs{X}}
 \newcommand{\BF}{\bs{F}}
 \newcommand{\BC}{\bs{C}}
 
 \newcommand{\BN}{\bs{N}}
 \newcommand{\BP}{\bs{P}}
 \newcommand{\BK}{\bs{K}}
 
 \newcommand{\Bb}{\bs{b}}
 \newcommand{\WITH}{\quad\text{with}\quad}
 \newcommand{\Btau}{{\boldsymbol{\tau}}}
 \newcommand{\Div}{\mbox{Div}}
 \newcommand{\Bzero}{{\boldsymbol{\mathit 0}}}
 \newcommand{\req }[1]{(\ref{#1})}
 
 \DeclareMathAlphabet{\Bgothic}{U}{euf}{b}{n}
 \DeclareRobustCommand{\BfrakU}{{\Bgothic U}}
 
 \DeclareRobustCommand{\BfrakP}{{\Bgothic P}}

 \newcommand{\fterm}[1]{\fbox{$\displaystyle#1$}}
 \newcommand{\noii}[1]{{\textcolor{black}{#1}}}
 
 \usepackage{soul} 
 \usepackage{dsfont} 
 \input{commands}

\begin{document}
 	
 	\title*{
Global-Local Forward Models within Bayesian Inversion for Large Strain Fracturing in Porous Media  
}
 	\titlerunning{Global-local within Bayesian Inversion} 
 	\author{Nima Noii, Thomas Wick, Amirreza Khodadadian}
 	\institute{Thomas Wick and Amirreza Khodadadian \at Leibniz University Hannover, Institute of Applied Mathematics and Cluster of Excellence PhoenixD, Welfengarten 1, 30167 Hannover,
 		Germany \email{khodadadian@ifam.uni-hannover.de}~~ \email{thomas.wick@ifam.uni-hannover.de}\\
 		Nima Noii  \at Leibniz University Hannover, Institute of Continuum Mechanics, An der Universität 1, 30823 Garbsen, Germany \email{noii@ikm.uni-hannover.de}
 		}
 	\maketitle
 	\vspace{-2.75cm}
 	\abstract{  		
 	In this work, Bayesian inversion with global-local forwards models is used to identify the parameters based on hydraulic fractures in porous media. 
It is well-known that using the Bayesian inversion \blue{to identify material parameters} is computationally expensive. \blue{Although each sampling may take more than one hour, thousands of samples are required to capture the target density.}
 	Thus, \blue{instead of using fine-scale high-fidelity simulations, we use a non-intrusive global-local (GL) approach for the forward model.}
We further extend prior work to a large deformation setting based on the Neo-Hookean strain energy function. The resulting framework is described in detail and substantiated with some numerical tests.}
 	\vspace{-0.5cm}

 	\section{Introduction}
 	\label{sec:1}
 	 	 	\vspace{-0.2cm}
 	Phase-field fracture models   are employed to capture failure and cracks in structures, alloys, and poroelastic media. The coupled model is based on solving the elasticity equation and an Allen-Cahn-type phase-field equation. In hydraulic fracture, a Darcy-type equation is solved to capture the pressure profile. Solving this coupled system of equations is computationally expensive. Indeed, to provide an accurate estimation (compared to the measurement) a very fine mesh profile is required. Of course, the time-dependent and nonlinear nature of the problem \blue{gives rise to more}  complexity. Another challenge is related to the computational, mechanical, and geomechanical material parameters. They have an essential effect on the simulations; however, many of them can not be estimated experimentally.
 	 
   In  \cite{noii2022bayesian}, we used the Bayesian inversion  to identify the parameters based on hydraulic fractures of porous media. A fracture response is realized through a phase-field equation  \cite{bourdin2000numerical} (based on the seminal work \cite{francfort1998revisiting}). But that work is limited to small deformations. In the current study, we extend \cite{noii2022bayesian} towards a \textit{large strain} formulation \cite{aldakheel2021global,MieheMauthe2015}. 
   
   In consequence, the main objective is to utilize non-intrusive 
global-local models \cite{gendre2009non} that are originally 
based on non-overlapping domain decomposition \cite{toselli2004domain}  
to significantly reduce the computational 
cost in Bayesian inversion. In extension to our prior work, we introduce an adoption of the hydraulic phase-field fracture formulation of a material that undergoes large deformation in poroelastic media. Finally, ensemble Kalman filters are employed \blue{for the proposal adaption  in Bayesian inversion to} identify the mechanical material parameters once the multiscale approach is used to solve the forward model.

\section{Framework for failure mechanics in hydraulic fracture}\label{section20}
\begin{figure}[!b]
	\centering
	{\includegraphics[clip,trim=7cm 4cm 1cm 5.5cm, width=10cm]{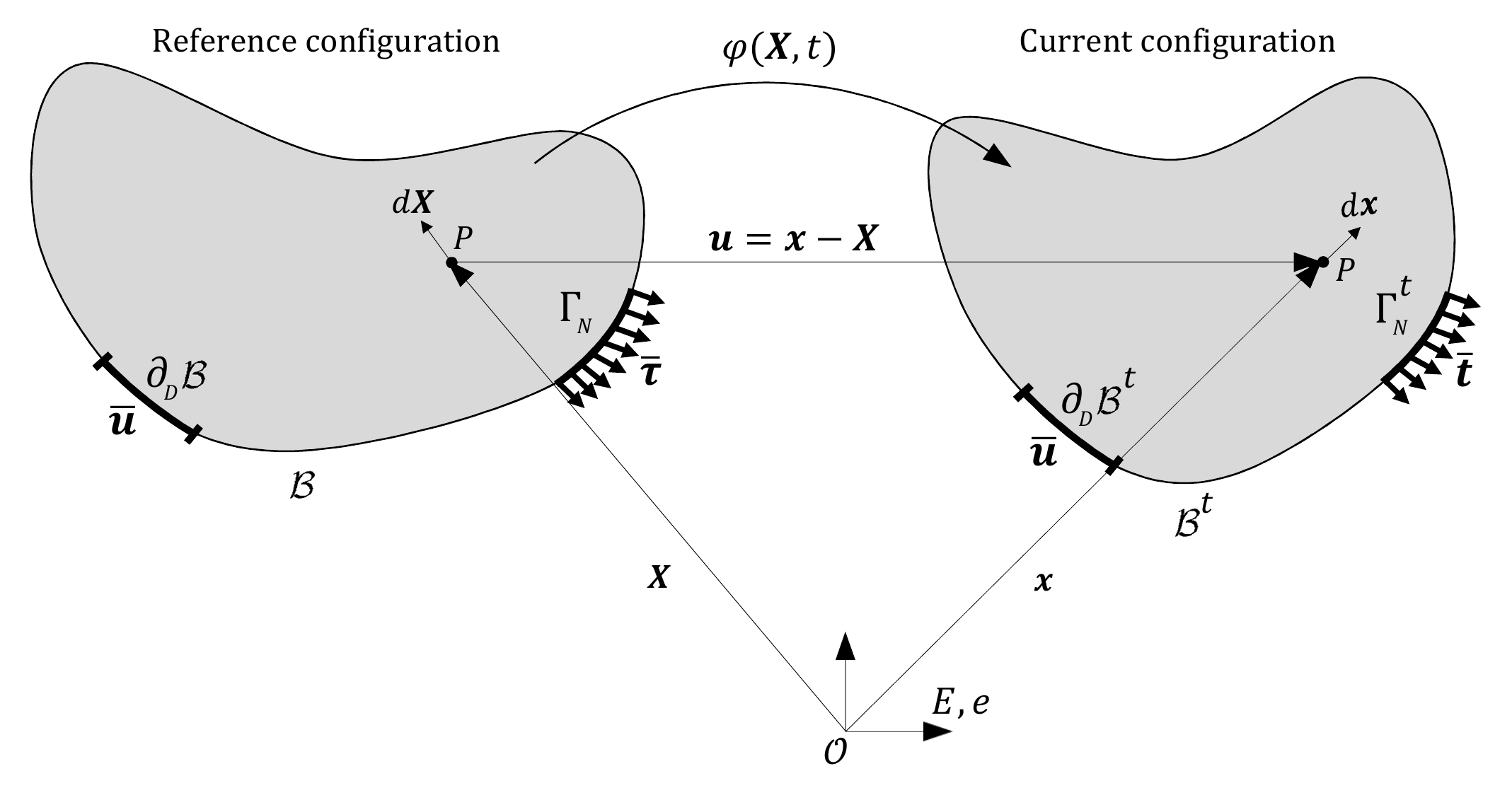}}  
	\caption{Setup of the notation for the configuration and motion  of the continuum body $\boldsymbol{\varphi}(\mathbf{X},t)$. The initial position $\BX$ in the undeformed configuration $\calB$ toward the current position $\Bx$ in the spatial configuration $\calB^t$ for the solid material undergoing finite strain.}
	\label{Fig_cm1}
\end{figure}

 	\vspace{-0.2cm}
 	Let us \blue{assume} $\calB\subset{{\mathbb{R}}}^{\delta}$ \blue{is the} solid computational domain (here $\delta=2$) with its surface boundary $\partial\calB$ and time $t\in \calT = [0,T]$. 
 	The given boundary-value problem (BVP) is a coupled multi-field system for the fluid-saturated porous media of the fracturing material. Since we are dealing in large strain setting, it is required to define the mapping between  the referential position $\BX$ towards spatial description $\Bx$ based on the motion $\Bvarphi$ of point $P$ at time $t$, see Figure \ref{Fig_cm1}. The media can be formulated based on a coupled three-field system. At material points $\Bx\in\mathcal{B}$ and time $t\in\calT$, the BVP solution indicates  the deformation field $\Bvarphi(\Bx,t)$ of the solid, the fluid pressure field $p(\Bx,t)$, and the phase-field fracture variable $d$ can be represented by
\begin{equation}
	\Bvarphi: 
	\left\{
	\begin{array}{ll}
		\calB \times \calT \rightarrow \calR^{\delta} \\
		(\BX, t)  \mapsto \Bx = \Bvarphi(\BX,t)
	\end{array}
	\right.~~~
	p:
	\left\{
	\begin{array}{ll}
		\calB \times \calT \rightarrow \calR \\
		(\BX, t)  \mapsto p(\BX,t)
	\end{array}
	\right.~~~
	d: 
	\left\{
	\begin{array}{ll}
		\calB \times \calT \rightarrow [0,1] \\
		(\BX, t)  \mapsto d(\BX,t)
	\end{array}
	\right .
	\label{phi-p-d-fields}
\end{equation}
 	
 	Here, $d(\Bx,t)=0$ and $d(\Bx,t)=1$ are referred to as the unfractured and completely fractured parts of the material, respectively. The coupled BVP is formulated through three specific primary fields to illustrate the hydro-poro-elasticity of fluid-saturated porous media by
 	\begin{equation}
 		\mbox{Global Primary Fields}: \BfrakU := \{ \Bvarphi, p, d \}.
 		\label{global-fields}
 	\end{equation}

\vspace{-0.85cm}
\subsection{Elastic contribution}
\vspace{-0.25cm}
The elastic density function is formulated through a Neo-Hookean strain energy function for a compressible isotropic elastic solid
\begin{equation}
	{W}_{elas}(\BF, d)= g(d)\; {\psi}_{elas}(\BF) 
	\WITH
	{\psi}_{elas}(\BF) = \frac{\mu}{2}\Big[(\BF:\BF -3) + \frac{2}{\beta} (J^{-\beta} -1) \Big]\;,
	\label{elas-part}
\end{equation} 
such that the shear modulus $\mu$ and the parameter $\beta:=\beta(\nu) = 2\nu/(1-2\nu)$ with the Poisson number $\nu<0.5$ are used.  \noii{Here, the material deformation gradient of the solid denoted by $\BF(\BX) :=
	\nabla\Bvarphi(\BX,t)=\text{Grad}\Bvarphi$ with the Jacobian $J \!:=\!\mbox{det}[\BF]  > 0$ augmented with the symmetric right Cauchy-Green tensor $\BC = \BF^T \BF$ is used; for details the reader is referred to \cite{aldakheel2021global,noii2022bayesian}. }
\noii{We note that the quadratic function ${g}(d) = (1-d)^2+\kappa$ is denoted as a degradation function, with $\kappa\approx10^{-8}$ that is chosen as a sufficiently small quantity.}
According to the classical Terzaghi theorem, the constitutive modeling results in the additive split of the stress tensor $\BP$ to effective mechanical contribution and fluid part as
\begin{equation}
	\begin{array}{ll}
		\BP(\BF, p, d) &:= \frac{\partial{	{W}_{elas}}}{\partial\BF }= g(d) \BP_{eff}(\BF) -
		B p J \BF^{-T}
		\WITH
		\BP_{eff} = \mu \big[ \BF - J^{-\beta} \BF^{-T}
		\big]\;\blue{.}
	\end{array}
	\label{p-piola-stresses}
\end{equation}
Here, the first Piola-Kirchoff stress tensor $\BP$ is derived from the first-order derivative of the pseudo-energy density function ${W}_{elas}$ given in \req{elas-part}. Thus, the balance of linear momentum for the multi-field system \noii{prescribed through body force $\overline\Bb$ reads}
\begin{equation}
	\fterm{
		\Div\,\BP(\BF, p, d) + \overline\Bb = \Bzero.
		\label{equil:defo}
	}
\end{equation}
\subsubsection{Fluid contribution}
The fluid volume flux vector \noii{$\BcalF$} is described through the negative direction of the
gradient of the fluid pressure $\nabla p$ and permeability based on Darcy-type fluid's
\begin{equation}
	\BcalF := - \BK(\BF,d)\; \nabla p.
\end{equation}
Here, the second-order permeability tensor \noii{$\BK(\BF,d)$}, following  \cite{MieheMauthe2015}, is additively decomposed into the permeability tensor into a Darcy-type flow for the unfractured porous medium ${\boldsymbol{K}_{Darcy}}$ and Poiseuille-type flow in a completely fractured material ${\boldsymbol{K}_{frac}}$ by
\begin{equation}
	\begin{array}{ll}
		\boldsymbol{K}(\BF,d) &=
		{\boldsymbol{K}_{Darcy}}(\BF)+d^{\zeta}
		{\boldsymbol{K}_{frac}}(\BF)
		\;, \\ [3mm]
		{\boldsymbol{K}_{Darcy}}(\BF) &= \frac{K}{\eta_F} J \BC^{-1}\;, \\ [3mm]
		{\boldsymbol{K}_{frac}}(\BF) &= K_c\; \omega^2\; J \big[
		\BC^{-1} - \BC^{-1} \BN \otimes  \BC^{-1} \BN
		\big].
	\end{array}
\end{equation}
 Here,  $K_D$ is the isotopic intrinsic permeability of the pore space, \noii{$K_c$ is the spatial permeability in fracture,} $\eta_F$ is the dynamic fluid viscosity, and $\zeta \ge 1$ is a
permeability transition exponent. Following \cite{MieheMauthe2015}, the so-called  {crack aperture (or the crack opening deformation) defined through $\omega = (\lambda_{\bot}-1) h_{e}$ in terms of the stretch orthogonal to the crack surface $\lambda^2_{\bot} ={\nabla d \cdot \nabla d}/{\nabla d \cdot \BC^{-1} \cdot \nabla d}$ and the characteristic element length $h_{e}$. Also, $\noii{\BN}=\nabla d/|\nabla d|$ 
	denotes the outward unit normal to the fracture surface, $h_{e}$ is the characteristic discretization size, and $\BI$ is an identity tensor. Thus, following  \cite{MieheMauthe2015,aldakheel2021global}, the fluid equation involve pressure files \blue{read}
\begin{equation}
	\fterm{
		\frac{\dot{p}}{M} + B \dot{J} - \bar{r}_F+ \Div[\BcalF] = 0\;\blue{.}
	}
	\label{pres-pde}
\end{equation}
\vspace{-0.75cm}
\subsubsection{Fracture contribution}
\vspace{-0.25cm}
The crack driving state function in the \blue{regularized} sense conjugate to crack phase-field denoted as ${D}(\Bvarphi,d,\Bx)$ for every point $\Bx$ in domain act as a driving force for the fracture evolution state reads
\begin{equation}
	\begin{aligned}
		{D}(\Bvarphi,d,\Bx):=\frac{2l}{G_c}(1-\kappa)\psi_{elas}(\BF).
		\label{eq2F}
	\end{aligned}
\end{equation}
\noii{Here, $G_c$ is the Griffith’s critical elastic energy release rate, and $l=2h_e$ is the regularization term.} Following \cite{miehe2015phase},  the local evolution of the crack phase-field equation in the given domain $\calB$ results in the third Euler-Lagrange differential system as
\begin{equation*}
	\fterm{
		(1-d)\calH-[d - l^2 \Delta d] = \eta\dot{d} \quad in~ \calB,
	}
	\tag{D}
\end{equation*}
augmented by the homogeneous Neumann boundary condition that is $\nabla d \cdot \Bn = 0$ on $\partial\calB$, with the maximum absolute value for the crack driving state $\calH = \max_{s\in [0,t]} D(\Bvarphi) \ge 0$ to avoid irreversibly. For different approach see \cite{noii2021quasi}. Thus, following our recent work \cite{noii2022bayesian}, the variational formulations for the three PDEs for the coupled poroelastic media of the fracturing material are 
\vspace{-0.1cm}
\begin{equation}
	\begin{array}{ll}
		\mathcal{E}_\varphi(\BfrakU, \delta \Bvarphi) &= \displaystyle\displaystyle \int_\calB \Big[ \BP:\nabla \delta \Bvarphi - \bar{\Bb} \cdot \delta \Bvarphi \Big] dV - \displaystyle \int_{\partial_N\calB} \bar{\Btau} \cdot \delta \Bvarphi  \; dA
		= 0 \ , \\ [4mm]
		\mathcal{E}_p(\BfrakU, \delta p) &= \displaystyle \int_\calB \Big[\Big(\frac{1}{M}(p-p_n) + B (J-J_n) -\Delta t \;\bar{r}_F
		\Big)\delta p + (\Delta t \;\BK \;\nabla p) \cdot \nabla \delta p \Big] dV \\ [3mm]
		& + \displaystyle \int_{\partial_N\calB} \bar{f} \;\delta p\; dA = 0 \ , \\ [4mm]
		\mathcal{E}_d(\BfrakU, \delta d) &= \displaystyle \int_\calB \Big[ \Big(2\psi_c \;d + 2(d-1) \calH \Big) \delta d  + 2\psi_c \;l^2 \; \nabla d \cdot \nabla \delta d
		\Big] dV = 0 \ .
	\end{array}
	\label{weakForm}
\end{equation}
This set of equation is now written in the abstract form through:  $\texttt{SS}(\BfrakU)$.
 	\begin{figure}[!t]
	\centering
	{\includegraphics[clip,trim=1cm 6cm 5cm 14cm, width=11cm]{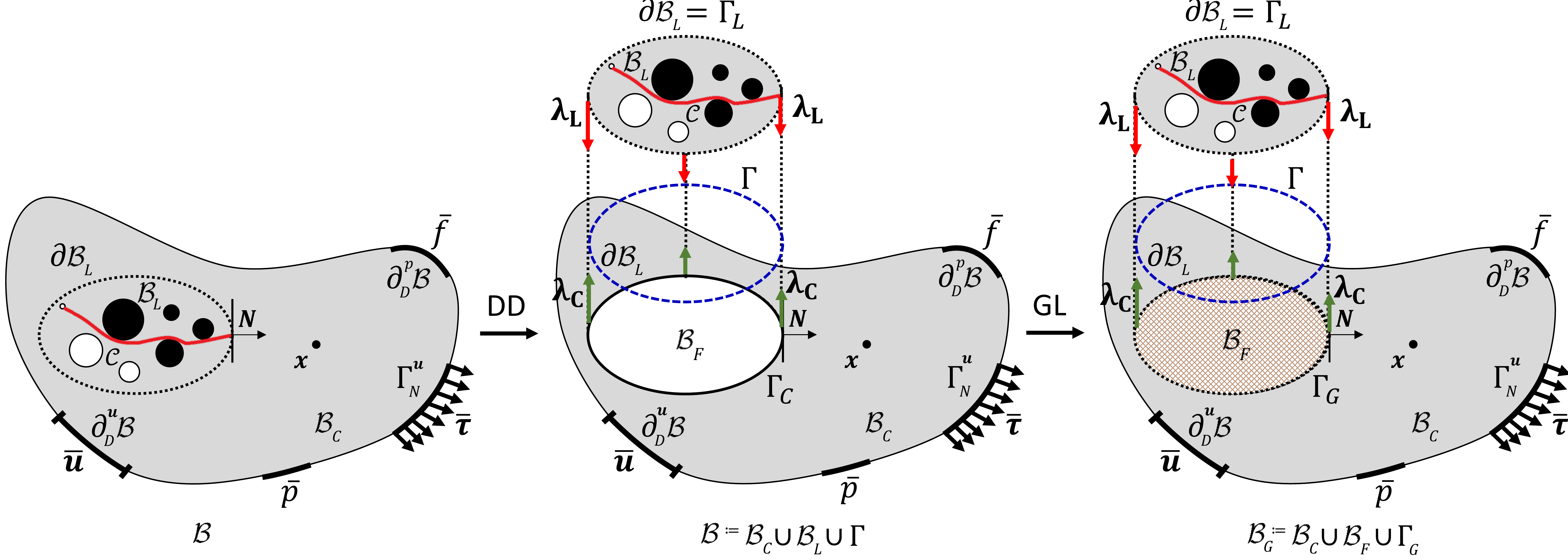}}  
	\caption{Configuration and loading setup of the single-scale BVP (left). Middle/right: global-local configuration, by the fictitious domain $\calB_F$ through filling the gap between $\calB_C$ and $\calB_L$ with a same constitutive modeling and discretization of $\calB_C$ such that its unification is a so-called global domain $\calB_G:=\calB_C\cup\Gamma_G\cup\calB_F$.
	}
	\label{Fig1}
\end{figure}
\vspace{-0.5cm}

\section{Multiscale modeling via a non-intrusive global-local method}\label{section2}
 	\vspace{-0.3cm}
 	The previously introduced system of equations 
for single-scale analysis in  \req{weakForm}
 	for the coupled problem of poroelasticity and fracture is further extended towards the global-local (GL) method now. Following \cite{aldakheel2021global,GeNoiiAllLo18}, the GL formulation is rooted in domain decomposition (e.g., \cite{toselli2004domain}) by distinguishing the original domain  into coarse and fine discretizations, see Figure \ref{Fig1}. 
To couple the domains, namely global and local domains, we have introduced an additional auxiliary interface denoted as $\Gamma$ between two disjoint domains in poroelastic media (see \cite{aldakheel2021global}),  and thus corresponding unknown fields, see Figure \ref{Fig1}.  These additional fields are the interface deformation $\Bvarphi_\Gamma(\Bx,t)$ and pressure $p_\Gamma(\Bx,t)$ on auxiliary interface and their corresponding traction forces $\{\Blambda_L^\Bvarphi, \Blambda^\Bvarphi_C\}$ and $\{\lambda_L^p, \lambda_C^p\}$ that are introduced as Lagrange multipliers. These results in a set of coupling equations at the interface by
	\begin{equation}
	\left\{
	\begin{tabular}{ll}
		${\Bvarphi}_L(\BX,t)  = {\Bvarphi}_{\Gamma}(\BX,t) $ & \mbox{at} ${\BX}\in{\Gamma}_{L}$, \\[0.1cm]
		${\Bvarphi}_G(\BX,t)  = {\Bvarphi}_{\Gamma}(\BX,t) $ & \mbox{at} ${\BX}\in{\Gamma}_{G}$, \\[0.1cm]
		$\Blambda_L^\varphi(\BX,t) + \Blambda^\varphi_C(\BX,t) = \Bzero$ & \mbox{at} ${\BX}\in{\Gamma}$,
	\end{tabular}
	\right.
	\ \mbox{and} \
	\left\{
	\begin{tabular}{ll}
		$p_L(\BX,t) = p_{\Gamma}(\BX,t)$ & \mbox{at} ${\BX}\in{\Gamma}_{L}$, \\[0.1cm]
		$p_G(\BX,t) = p_{\Gamma}(\BX,t)$ & \mbox{at} ${\BX}\in{\Gamma}_{G}$, \\[0.1cm]
		$\lambda_L^p(\BX,t) + \lambda^p_C(\BX,t) = \textit{0}$ & \mbox{at} ${\BX}\in{\Gamma}$.
	\end{tabular}
	\right.
\end{equation}
 Now, the multi-physics problem for the global-local approach is described through  eleven primary fields to characterize the hydro-poro-elasticity of fluid-saturated porous media at finite strains by
 \begin{equation}
 	\mbox{Extended Primary Fields}: \BfrakP := \{ \Bvarphi_G, \Bvarphi_L, p_G, p_L, d_L, \Blambda^\Bu_C, \Blambda^\Bu_L, \lambda^p_C, \lambda^p_L, \Bu_\Gamma, p_\Gamma \}
 	\label{gl-fields}
 	\ .
 \end{equation}
 Herein, a global constitutive model behaves as a poroelastic response, abbreviated as E(elastic)-P(pressure), which is augmented with a \textit{single local domain} and behaves as a poroelastic material with fracture response, abbreviated as E(elastic)-P(pressure)-D(damage). 
The resulting final algorithm is based on our prior work \cite{aldakheel2021global,noii2022bayesian}.
 
 \vspace{-0.7cm}

\section{Bayesian inversion for parameter estimation }\label{section3}
 	 \vspace{-0.3cm}
 	  In this study, we have use MCMC (Markov chain Monte Carlo) techniques to identify the material parameters in the hydraulic porous medium phase-field fracture setting.
The latter is solved with the previously described GL approach. 
In general, we can employ the following probabilistic model to update the available prior information according to the forward model (here considers the phase-filed fracture) and a reference observation (arising from measurement, or a synthetic observation). First, we introduce the following statistical model
 	     	 	\vspace{-0.2cm}
 	 \begin{align}
 	 	\label{BE}
 	 	\mathbb{M}=\mathcal{P}(\xx,\chi)+\varepsilon.
 	 \end{align}
 	 Here $ \mathbb{M}$ refers to the reference observation arising from the experimental data (a measured value) and $\mathcal{P}$ considers to the model response related to $\chi$ \blue{a set of $d$-dimensional material parameters}. Furthermore, $\xx\in \mathbb{R}^\delta$ and $\varepsilon$ indicates the measurement error. It assumed to have Gaussian independent and identically distributed error $\varepsilon\sim\mathcal{N}(0,\sigma^2\,I)$, having the parameter $\sigma^2$. 
 	 Since $\mathcal{P}$ in \req{BE} is a model response which results in our computation, such that in our presented model can be approximated through
 	 \begin{align*}
 	 	\texttt{signle-scale:} \quad \mathcal{P}\approx\mathcal{P}^{\text{SS}} 
 	 	\quad \text{or} \quad
 	 	\texttt{global-local:} \quad \mathcal{P}\approx\mathcal{P}^{\text{GL}} , 
 	 \end{align*} 
 	 corresponds to equations (SS) and (GL), respectively. Thus, \req{BE} becomes as
 	 \begin{align}
 	 	\mathcal{M}=\mathcal{P}^\bullet(\Theta)+\varepsilon,\WITH \bullet\in\{\text{SS,GL\}}.
 	 \end{align} 
  
  	 \begin{figure}[!t]
 	\centering
 	\includegraphics[width=0.9\textwidth]{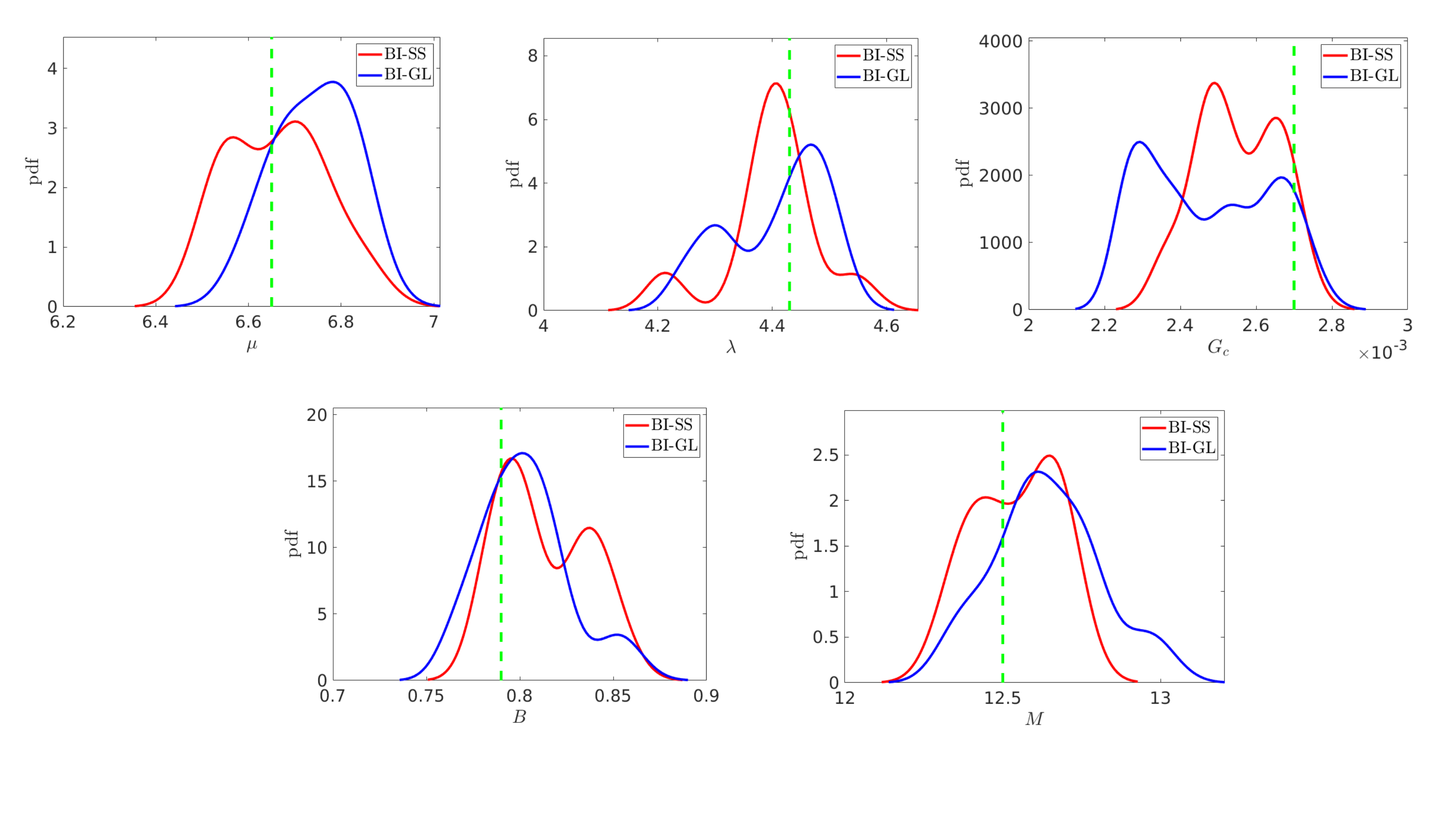}
 	\vspace{-0.8cm}
 	\caption{The pdf of posterior density of the material parameters using the BI-GL and BI-SS approaches for  fracture. The true values are shown with a dashed green line.}
 \end{figure} 	
	\begin{figure}[!b]
	\centering
	{\includegraphics[clip,trim=0cm 13.5cm 10.5cm 14cm, width=9cm]{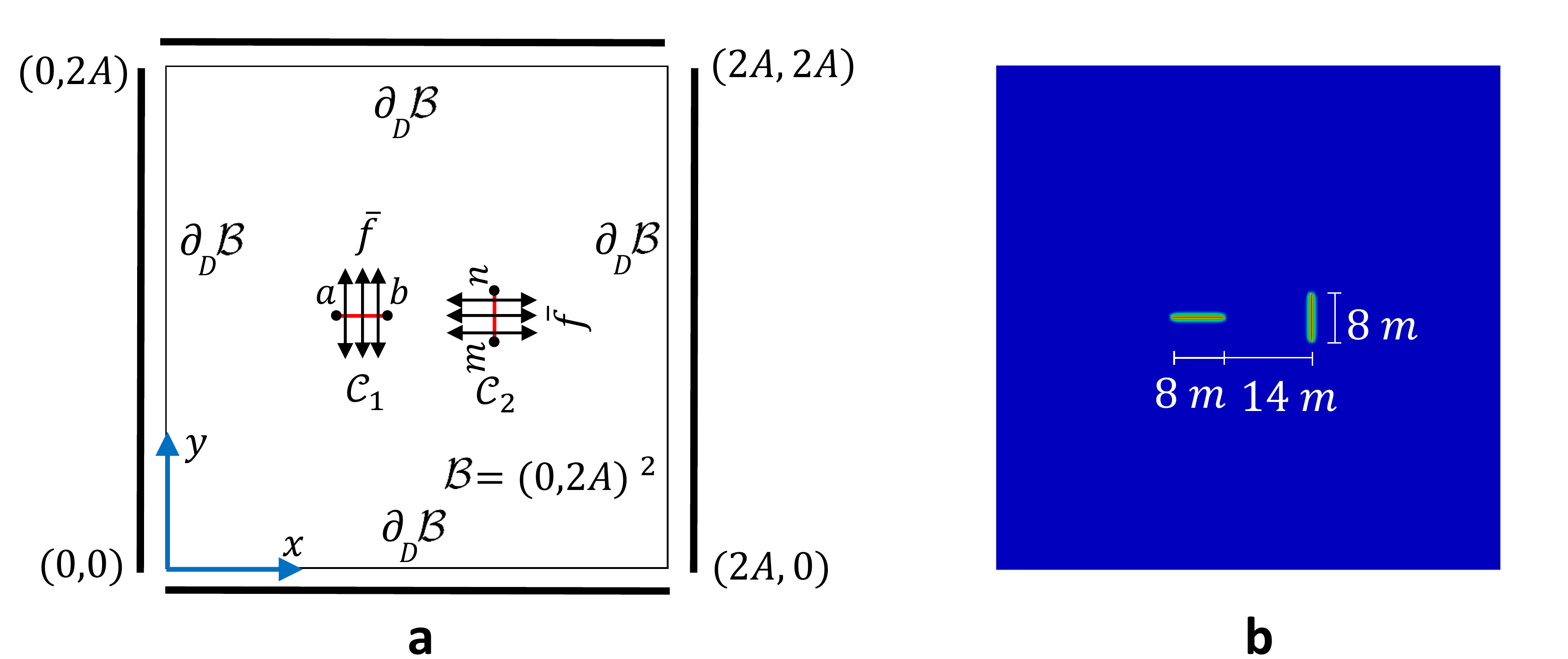}}  
	\caption{Joining of two cracks driven by fluid volume injection. (a) Geometry and boundary conditions; and (b) described crack phase-field $d$ as a Dirichlet boundary conditions at $t=0~s$.
	}
	\label{example2-a}
\end{figure}
	
 	 Despite the simplicity of the \blue{Metropolis-Hastings} algorithm, it is not suitable for complicated cases, specifically when several parameters should be estimated (multi-dimensional domains). In this study, we have used MCMC with ensemble-Kalman filter, see for a detailed discussion \cite{noii2022bayesian}. Ensemble Kalman filter (EnKF) indicates the error covariance matrix by a large random ensemble of model observations.
 	 \blue{Here, to achieve a reliable estimation of posterior density, a Kalman gain is computed using the mean and the covariance of the prior density and the cross-covariance between material parameters and observations. Using an ensemble-Kalman filter, we adopt the proposal density with $\chi^\star=\chi^{j-1}+\Delta \chi$, where  $\Delta \theta$ is the jump of Kalman-inspired proposal. Afterwards, we update the candidate via $\Delta \chi=\mathcal{K}\left( y^{j-1}+s^{j-1}\right)$. The Kalman gain is computed by $\mathcal{K}=\C_{\theta M}\left(\C_{MM}+\mathcal{R}\right)^{-1}$, where $\C_{\theta M}$ is the covariance matrix between the unknowns and the model response, $\C_{MM}$ denotes the covariance matrix of the PDE-based model, and $\mathcal{R}$ is the measurement noise covariance matrix \cite{zhang2020improving}. Moreover, $y^{j-1}$ is the residual of candidates w.r.t the model and  $s^{j-1}\sim\N(0,\mathcal{R})$  relates to the density of measurement.  Denoting $\texttt{obs}$ as an observation, $y^{j-1}=\texttt{obs}-f(\theta^{j-1})$. We refer the reader to \cite{noii2022bayesian_B} for more details and the codes.}

   Thus, we are now able to use Bayesian inversion to identify the fracking process using multiscale approach material parameters that cannot be measured with usual techniques.

 	\vspace{-0.75cm}
 
\section{Numerical example}
\vspace{-0.25cm}
In this section, we investigate a numerical test with the main goal that Bayesian inversion yields accurate parameter identifications at a cheap cost of the governing global-local phase-field solver. The mechanical and geomechanical descriptoion of the parameters is given in \cite{noii2021bayesian}.
 	In the following, a BVP is applied to the square plate shown in Figure \ref{example2-a}. The geometry and boundary conditions are from \cite{aldakheel2021global}. The \blue{single-scale (SS) model results} considering the phase-field and pressure are given in Figure \ref{SS}. Then, we employ our global-local approach, 
with findings shown in Figure \ref{GL}. Figure \ref{both} shows the load-displacement curve for both approachs, indicating the accuracy of the GL approach. Finally, the computational costs of both approaches using the Bayesian setting is given in Table \ref{table1}, denoting the significant efficiency of the domain decomposition technique.

 	 	\begin{figure}
 		\centering
 		\includegraphics[width=0.9\textwidth]{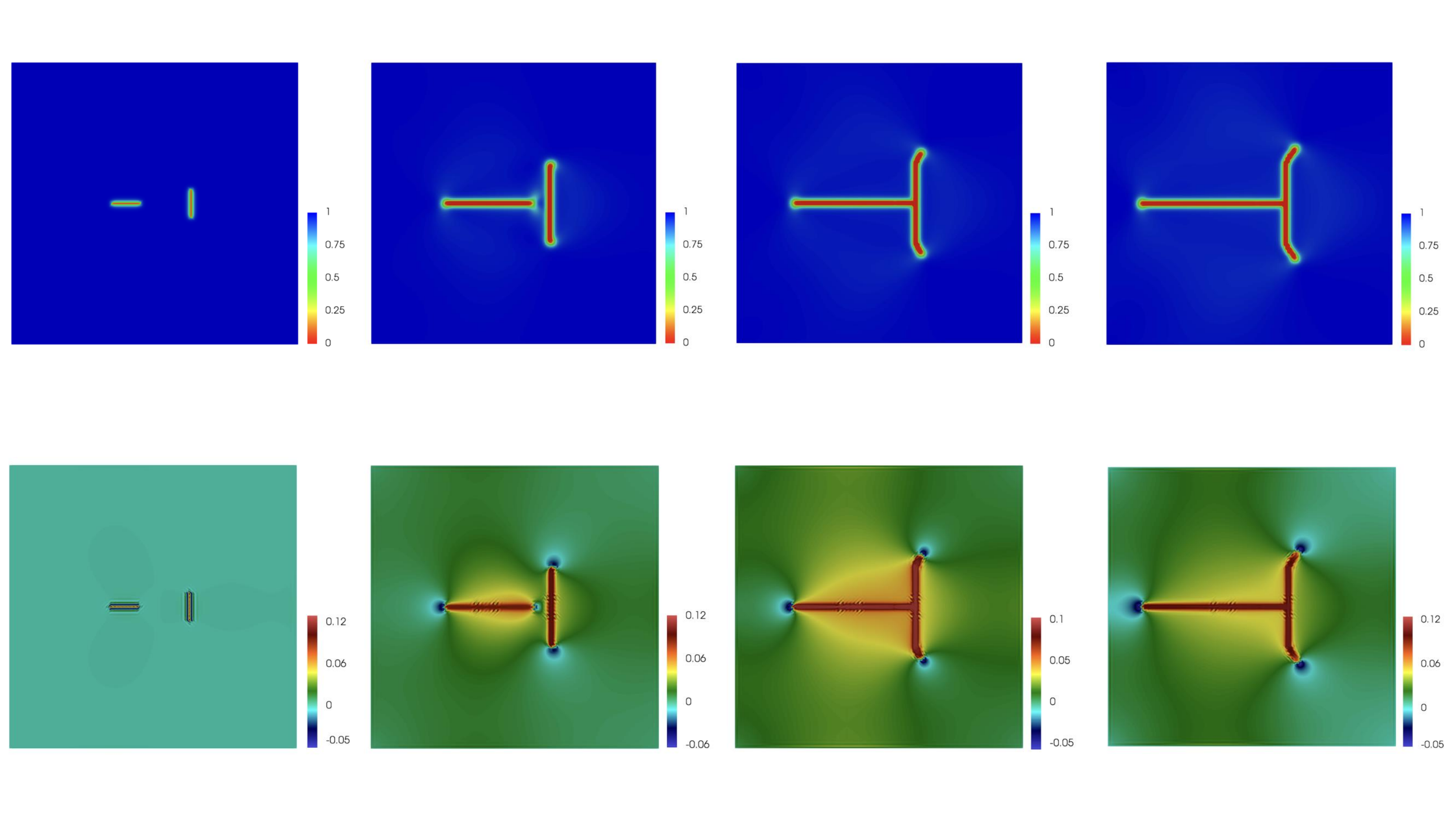}
 		\vspace{-0.5cm}
 		\caption{The evolution of the phase-field (first line) and pressure (second line) for different fluid injection time, i.e., $t\in\,[0.1, 10, 15, 20]$ seconds using SS model.
 		}
 	\label{SS}
 	\end{figure}

 		\begin{figure}
 		\centering
 		\includegraphics[width=0.9\textwidth]{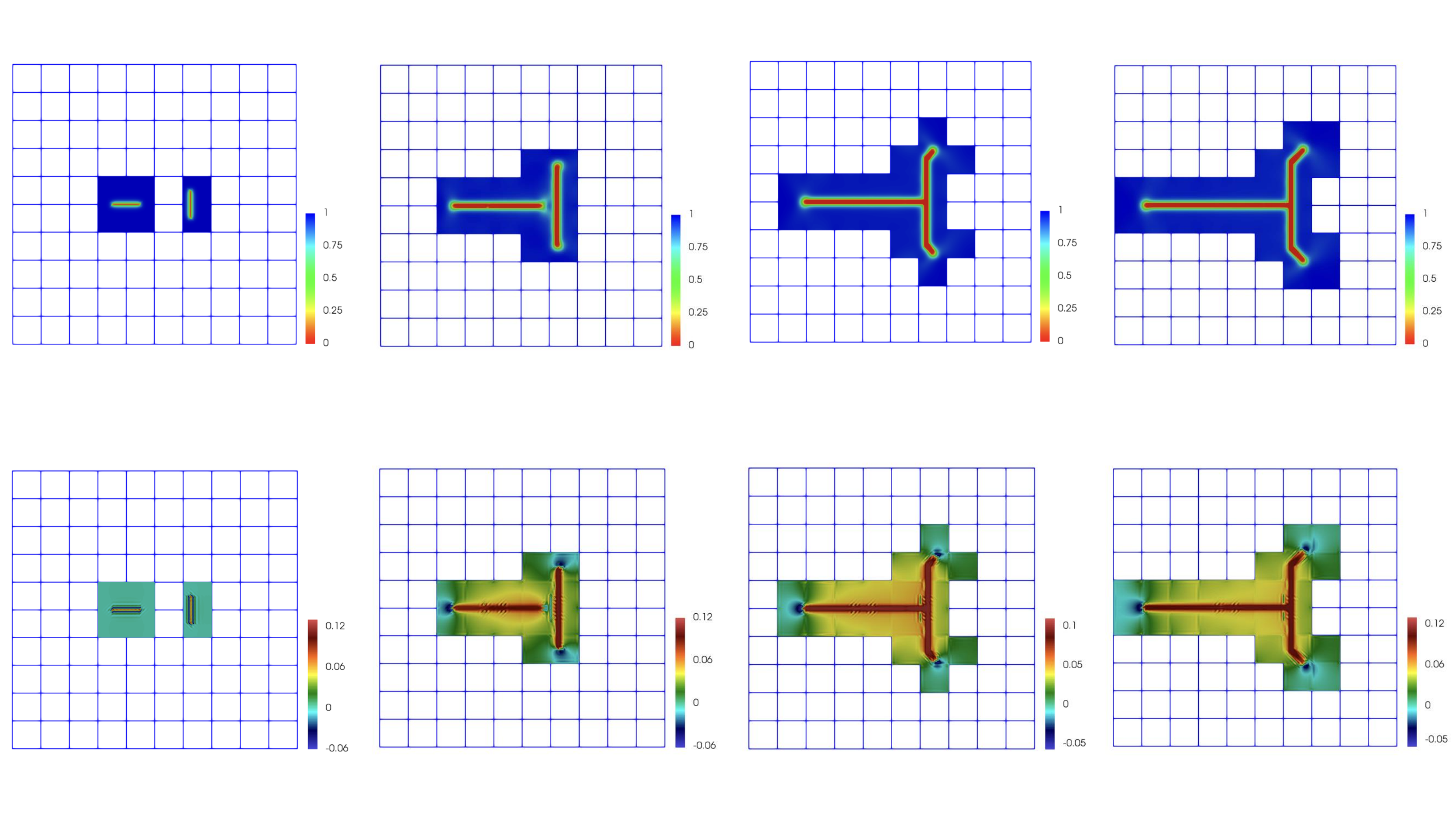}
 		\vspace{-0.5cm}
 		\caption{The evolution of the phase-field (first line) and pressure (second line) for different fluid injection time, i.e., $t\in\,[0.1, 10, 15, 20]$ seconds using GL model.}
 		\label{GL}
 	\end{figure}


 	\begin{figure}
 		\centering
 		\includegraphics[width=0.9\textwidth]{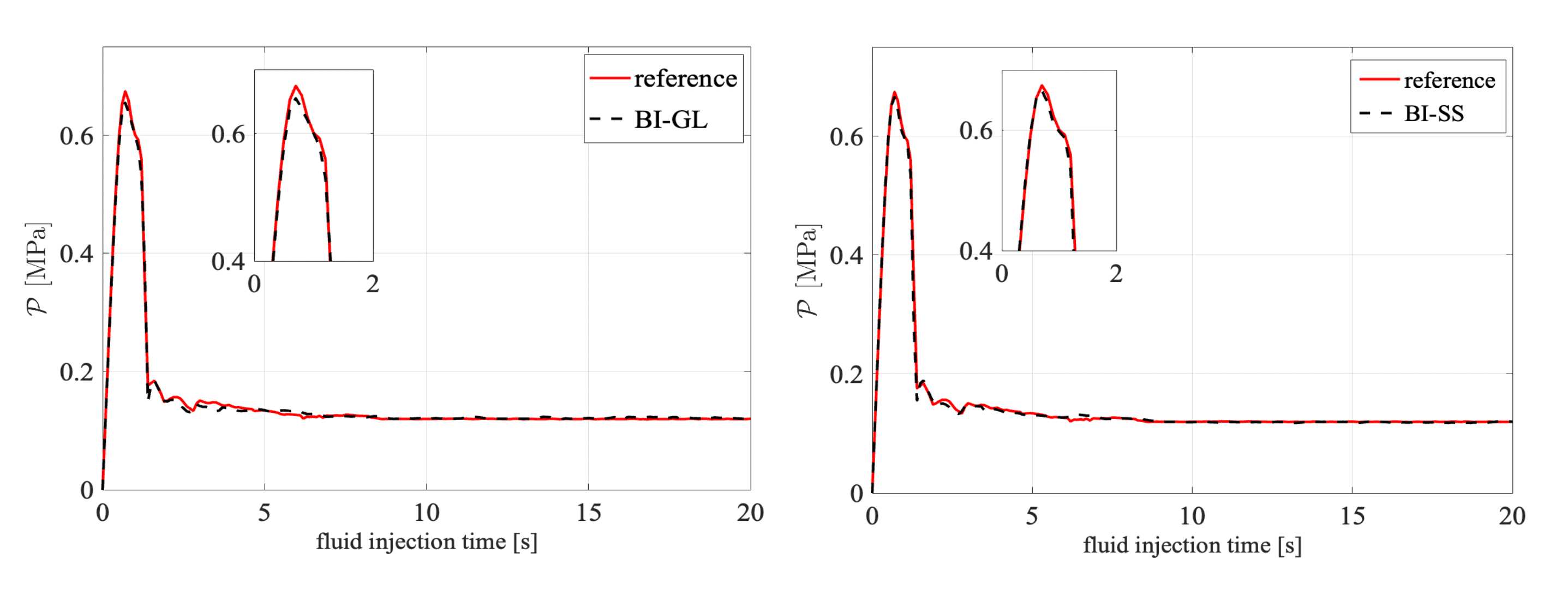}
 		 		 		\vspace{-0.4cm}
 		\caption{A comparison between the maximum pressure obtained by the true values (the reference observation) and the mean value of posterior density of BI-GL (left) and BI-SS (right).}
 		\label{both}
 	\end{figure}


	\begin{table}\label{time2}
 	\small
 	\caption{A comparison between the computational costs of BI-SS and BI-GL approaches for hydraulic fracture. The unit is given in seconds.} 
 	\centering
 	\begin{tabular}{|lccccc|}
 		\hline
 		Model&  \quad min $T$\quad        &  \quad max $T$\quad  &   \quad mean $T$\quad      &   \quad $\displaystyle \sum\; T$ \quad&\quad ratio $T$ \quad       \\[0.5mm]\hline
 		$\mathrm{BI-SS}$    &\quad   5\,645 &  5\,767 &  5\,704 & 1.14$\times 10^{6}$ &19.47\\[4mm]
 		$\mathrm{BI-GL}$ & \quad 277  & 296.2  & 287.1& 5.75$\times 10^{4}$  & -- \\[4mm]
 		\hline
 	\end{tabular}
 \label{table1}
 \end{table}

 \vspace{-0.6cm}
 	\section{Conclusion}
 	 \vspace{-0.2cm}
 	\label{conclusions}
 	In this study, we extended a global-local (GL) approach for phase-field fracture as the PDE-based model with Bayesian inversion. We applied the proposed idea to hydraulic fracturing within poromechanics concept, for material undergoing large deformation. 
For our numerical example, Bayesian inversion
 	using GL is $20$ times faster than the signle-scale model, while the
 	accuracy is similar.

\vspace{-0.25cm}
 	\begin{acknowledgement}
 		N. Noii  acknowledges the Priority Program Deutsche Forschungsgemeinschaft \texttt{DFG-SPP 2020} within its second funding phase. T. Wick and A. Khodadadian acknowledge  
 		the DFG under Germany Excellence Strategy within  the Cluster of Excellence PhoenixD (EXC 2122, Project ID 390833453.
 	\end{acknowledgement}

 \vspace{-1.3cm}
 	\bibliographystyle{abbrv}
 	\bibliography{lit}
 	
 \end{document}

%% file: commands.tex
\definecolor{lightblue}{RGB}{0, 142, 255}

\definecolor{sblue}{RGB}{0, 10, 255}